\documentclass[12pt,a4paper]{article}
\usepackage[T1]{fontenc}
\usepackage{lmodern}
\usepackage[latin1]{inputenc}
\usepackage{latexsym}                 
\usepackage{color}                 
\usepackage{amssymb}
\usepackage{amsmath}
\usepackage{amsthm}
\usepackage{graphicx}
\usepackage{hyperref}
\usepackage{verbatim}
\usepackage{mathptmx}      
\def\eref#1{(\ref{#1}%
)}
\def\RSref#1{\ref{#1}%
}
\def\RSlabel#1{\label{#1}%
}
\def\RScite#1{\cite{#1}%
}

\newcommand{\bql}[1]{%
\begin{equation}\label{#1}%
}

\def\filename#1{}

\newcommand{\eq}{\end{equation}}
\def\fa{\hbox{ for all }}

\def\b1{\mathbf 1}

\newcommand{\R}{\ensuremath{\mathbb{R}}}

\def\biglf{\par\bigskip\noindent}
\newtheorem{theorem}{Theorem}
\newtheorem{corollary}{Corollary} 

\usepackage{amsmath}
\usepackage{upgreek}
\def\biglf{\par\bigskip\noindent}
\begin{document}
\begin{center}
  {\large \bf An Approximation Theorist's View on\\~\\
    Solving Operator Equations}\\~\\
    -- with special attention to Trefftz, MFS, MPS, and DRM methods --\\~
\biglf 
Robert Schaback\footnote{%
Institut für Numerische und Angewandte Mathematik,
Universit\"at G\"ottingen, Lotzestra\ss{}e 16--18, 37083 G\"ottingen, Germany,
{\tt schaback@math.uni-goettingen.de}}
\\
\end{center}
\par\noindent
{\bf Abstract}:
When an Approximation Theorist looks at well-posed
PDE problems or operator equations,
and standard solution algorithms like Finite Elements, Rayleigh-Ritz
or Trefftz techniques, methods of fundamental or particular solutions and
their combinations, they boil down to approximation problems
and stability issues. These two can be handled by Approximation Theory,
and this paper shows how, with special applications
to the aforementioned algorithms. The intention is that
the Approximation Theorist's viewpoint is helpful fur readers
who are somewhat away from that subject.
\section{Introduction}\RSlabel{SecIntro}

Whenever a specific unknown
    function $u^*$ is to be numerically constructed
    from whatever known information $D(u^*)$ about it, an
    Approximation Theorist will first look at {\em trial spaces}
    that can approximate the function well, including its data.
    Whatever the potential numerical recipes 
    are, the resulting errors for the calculated trial functions
    $\tilde u$ should always be
    comparable to the achievable error when approximating the
    true solution directly from the trial space, because that
    error cannot be improved.
    \biglf
    To make this argument work, operator equations and numerical
    algorithms come in the way, unfortunately,  and need to be surpassed,
    but they are not directly relevant to the Approximation Theorist's
    argument. Consequently, this
    viewpoint does not really care for the PDE or operator equation
    problem, and rightfully so,
    because it turns out that
    convergence rates are not PDE-dependent, if {\em stability}
    can be guaranteed. Then they depend only on
the obtainable approximation error, and the latter depends on the true
solution, its smoothness, and the chosen trial space.
The PDE problem does not matter, as long as
stability prevails, in a sense to be worked out.
\biglf
Consequently, an Approximation Theorist
will use solution techniques that are approximation problems as well,
and then the error analysis turns out to be quite simple, as will be shown.
The basic obstacles are stability problems that will be handled
via well--posedness of the operator equations, but then
convergence rates can be played back to Approximation Theory, as expected.
\biglf
This paper approaches these goals step by step, summarizing
results from \RScite{schaback:2010-2,schaback:2015-4,schaback:2015-3},
in a more concise form than before,
and applying them later explicitly
to Trefftz, MFS, MPS, and DRM
cases in specialized sections. Readers may be puzzled by the fact that
Approximation Theorists avoid linear systems of equations and
prefer to work in terms of spaces, not bases. But it will be clear why.
It simplifies things and avoids additional instabilities.
\biglf
The first step concerns {\em analytic theory}, starting from  a
general formulation
of {\em well-posed}
linear operator equations that allow FEM, Trefftz, MFS, and MPS
methods. It is assumed that the analytic problem has a solution, and
that the operators themselves
are not discretized, i.e. differential and boundary operators are
not replaced by finite differences, but applied directly to trial functions.
The differential and boundary operators that are defining a PDE problem
are merged into one single {\em data map} $D\;:\;U\to F$
that maps a function $u\in U$ to
the values $D(u)\in F$ of the operators in question. Solving the
problem then  consists of inverting $D$. A typical case would be
$D(u)=(-\Delta u, u_{|_{\Gamma}})$ for a Poisson problem with Dirichlet data.
It will be assumed
that the problem is stated in the form $D(u)=f$
for a given $f\in F$, but it is also assumed to be
solvable by a function $u^*$, i.e. $f=D(u^*)$ holds
for some $u^*\in U$.
Furthermore, the problem should be well-posed in the sense
\bql{eqWPOp}  
\|u\|_{WP}\leq C_{WP}\|D(u)\|_F \fa u\in U,
\eq
i.e. each function should be continuously recoverable from its data.
The above {\em well-posedness norm} $\|.\|_{WP}$ should be weaker than the norm
$\|.\|_U$ on $U$ and 
is of central importance to the error analysis to follow.
In case of elliptic boundary value problems
with the Maximum Principle, the {\em well-posedness norm} will be the sup norm.
This finishes the PDE side. The rest of the argumentation
will not depend on the PDE anymore, 
once the data map, the spaces, their norms, and the true solution are fixed
together with a well--posedness inequality \eref{eqWPOp}.
We shall illustrate this in Section \RSref{SecPDE}.
\biglf
The second step is {\em Approximation Theory}.
One should choose a finite-dimensional
{\em trial space} $U_r\subset U$ that 
is able to approximate the true solution $u^*$ well. If it does so, it will
also approximate the data $D(u^*)=f$ of $u^*$ well, i.e. one rather considers
the approximation of $f=D(u^*)$ by functions $f_r \in F_r:=D(U_r)\subset F$
and expects that Approximation Theory has good news about the
obtainable minimal error
\bql{eqVapp0}
\displaystyle{\inf_{f_r\in F_r}\|f-f_r\|_F=\|f-f_r^*\|_F=:\eta(f,F_r,F)}.
\eq
These approximations will hopefully 
determine the convergence rate of the algorithms that are to be defined
for PDE solving as well,  
and the rate will hopefully not depend on the PDE problem. 
The well-known standard example is that the classical
FEM convergence rate is the 
PDE-independent
convergence rate of piecewise linear approximations  in Sobolev spaces.
But this basic PDE-independence of convergence rates holds in general
and comes from Approximation Theory.
We add details in Section \RSref{SecAT}.
\biglf
The third step is {\em Theoretical Numerical Analysis}. If one solves
approximation problems 
instead of solving linear systems, namely by minimizing the residuals
$f-D(u_r)=D(u^*)-D(u_r)$ over all $D(u_r)$ 
in $F$, the convergence rates of the
previous Approximation Theory
step carry over to the numerical solution of the PDE problem,
and stability is automatically guaranteed by well-posedness.
This is easy to see via $F_r:=D(U_r)$ and
\bql{eqapperrbnd}
  \|u^*-u_r^*\|_{WP}
  \leq  C_{WP}\|f-D(u_r^*)\|_F\leq C_{WP}\;\eta(f,F_r,F),
\eq
and is quite satisfactory, as far as error analysis and convergence rates are
concerned. Stability problems do not arise as long as the analysis uses spaces,
not bases. Summarizing:
\begin{theorem}\RSlabel{TheConvContApp}
  Let a problem in the form $D(u)=f$ with a data map $D\;:\;U\to F$
  be given and assume it is
  well--posed in the sense
  of \eref{eqWPOp}. Pick a trial space $U_r\subset U$,
  form the space $F_r=D(U_r)$ and approximate $f$
  from $F_r$ in the norm on $V$ via \eref{eqVapp0}.
  Then the optimal solution $f_r^*=D(u_r^*)$ satisfies the error bound
  \eref{eqapperrbnd}, i.e. the error in the solution,
  measured in the well--posedness norm, is 
  is proportional to the error of approximating the
  data of the true solution. \qed 
\end{theorem} 
But the previous step employed
approximation in {\em function spaces}, and this is not easy to handle
in practice. It concerned a much too theoretical instance of Numerical
Analysis. 
Therefore we have to deal with
a fourth step, namely the problem of
{\em Discretization in Approximation Theory}, replacing functions by
finitely many values.
This has nothing to do with PDE theory again,
but it arises in the background of 
numerical methods like MFS, MPS, Trefftz, or DRM.
It is the part of Numerical Analysis that handles approximation problems
in function spaces and breaks them down to some form of Linear Algebra
or Optimization. Here, stability issues creep in through the back door.
It will be proven in Section \RSref{SecStab}
that certain approximation problems allow uniformly stable
discretizations, if functions are replaced by sufficiently many values,
and this applies to certain well-posed PDE problems
or operator equations in strong formulation.
The final sections illustrate
these results for  Trefftz
techniques, the Methods of Fundamental or Particular
Solutions, and the Dual Reciprocity Method. The experimental
paper \RScite{schaback:2008-7} has many numerical examples that
support an Approximation Theorist's view on these methods.  
\biglf
But before we go on, here is a seemingly trivial
practical consequence of Theorem
\RSref{TheConvContApp} concerning an a-posteriori error analysis:
\begin{corollary}\RSlabel{TheAPostErr}
For a well--posed problem in the above sense, assume
that a function $\tilde u\in U$  is produced by whatsoever method,
and assume that the norm $\|f-D(\tilde{u})\|_F$ can be calculated. Then
$$
  \|u^*-\tilde u\|_{WP}
  \leq  C_{WP}\|f-D(\tilde u)\|_F
  $$
is an error bound involving only the well-posedness constant. \qed
\end{corollary}
Unfortunately, the literature on operator equations only rarely
yields explicit upper bounds for $C_{WP}$. This topic deserves much more
attention in mathematical research.
\biglf
But the above argument allows
a fair comparison of different
numerical methods that produce numerical solutions
$\tilde u$ for the same well--posed problem. Even if $C_{WP}$
is not known, it is independent of the numerical techniques,
and the actual error $\|f-D(\tilde{u})\|_F$ of approximating the data
is a valuable information for the comparison of methods. Far too many
numerical papers insist on knowing the true solution, and produce examples with
unrealistically smooth true solutions. Instead, it suffices
to reproduce the data in the norm $\|.\|_F$
well, and hopefully better
than competing methods. If the data error
$\|f-D(\tilde{u})\|_F$ is presented in a paper, and if the error does not meet
the expectations of Approximation Theory, there is a serious
stability flaw in the presented method that needs special attention. 
\section{Operator Equations}\RSlabel{SecPDE}
We first specify which operator equations we shall consider,
and how standard PDE problems are subsumed.
\biglf
By Section \RSref{SecIntro}
an Approximation Theorist sees solving operator equations
as a numerically motivated
detour from the central problem of approximating the true solution $u^*$ by
functions from the trial space. The detour is necessary, because
one has only indirect information about the solution, e.g.
values of derivatives at certain places, or values of integrals of
derivatives against certain
test functions. These are the available {\em data} $D(u^*)$ of $u^*$.
\biglf
Thus an operator equation takes the form
\bql{eqDuf}
D(u)=f
\eq
for a {\em data map} $D\;:\;U\to F$  between
Banach spaces that is to be inverted. In particular,
we assume that $f=D(u^*)$ is given, i.e. the
problem is exactly solvable by some $u^*\in U$. The data map simply describes
what is known about the solution, e.g. the pair
$D(u):=((-\Delta u)_{|_{\Omega}}, u_{|_{\partial\Omega}})$
for a standard strong Poisson problem on a bounded domain $\Omega$.
\biglf
Besides solvability, we require {\em well-posedness} of the operator equation in
the sense of a {\em well-posedness inequality} \eref{eqWPOp}. This
implies continuous invertibility of $D$, but needs some
explanation, because it concerns the norms in $U$ and $F$.
The norm in $F$ also arises in the approximation problem
\eref{eqVapp0} and should not be too exotic. A typical bad case
arises when setting up PDE problems in Hölder or Sobolev space, because these
carry norms that are not easy to access numerically. 
\biglf
For classical strong Dirichlet problems for uniformly elliptic
self-adjoint second-order differential operators $L$
on compact domains $\Omega$ with boundary $\Gamma$,
there is a well-posedness inequality
\cite[p.14]{braess:2001-1}
\bql{eqinfWP}
\|u\|_{\infty,\Omega}
\leq \|u\|_{\infty, \Gamma} +C\|Lu\|_{\infty,\Omega}
\fa u\in U:=C^2(\Omega)\times C(\Gamma)
\eq
in the sup norm on $U$. Then we can choose the right-hand side as our norm in
$U$ and get well-posedness also in the norm on $U$.  The data space is
$F=C(\Omega)\times C(\Gamma)$ and carries manageable norms,
the data map being defined via
$$
D(u):=(Lu_{|_\Omega},u_{|_\Gamma}) \fa u\in U.
$$
Weak problems are different, because they have other data maps.
Authors should always consider strong and weak ``formulations''
as completely separate problems, not just two aspects of the same thing.
The difference comes up when we write the data maps in terms of infinitely many
conditions. The strong Poisson problem on $\Omega$
has infinitely many equations
\bql{eqPoiStrong}
\begin{array}{rclcrl}
-\Delta u(x) &=& f(x) &=&  -\Delta u^*(x),&x\in \Omega,\\
u(y) &=& g(y) &=&  u^*(y),&y\in \partial{\Omega}
\end{array}
\eq
while the corresponding weak problem consists of
\bql{eqPoiWeak}
\begin{array}{rccccl}
  (\nabla u,\nabla v)_{L_2(\Omega)} &=&  (v,f)_{L_2(\Omega)} &=&
  (\nabla u^*,\nabla v)_{L_2(\Omega)} ,&v\in H_0^1(\Omega),\\
u(y) &=& g(y) &=&  u^*(y),&y\in \partial{\Omega},
\end{array}
\eq
just another set of infinitely many equations. Note that by mixture
of the above cases one can pose very many different problems, with all kinds of
differential and boundary operators. But we refer the reader to
\RScite{schaback:2015-4}
for details on handling classical local and global weak problems
with this approach. 
\biglf
For Trefftz methods
\RScite{qin:2000-1, li-et-al:2008-1,kolodziej-zielinski:2009-1, li-et-al:2010-1}, one ideally has a homogeneous differential
equation and poses only boundary conditions. Then the data map $D$ can
consist only of the boundary values, and the space $U$ should be restricted
beforehand to homogeneous solutions.
The Method of Fundamental Solutions (MFS,
\RScite{mathon-et-al:1977-1,Bogomolny:1985-1,katsurada:1990-1,
  katsurada:1994-1,
  fairweather-et-al:1998-1, golberg-et-al:1998-2,
poullikkas-et-al:1998-1,chen-et-al:2008-1}) is a special case.
Details will be in Section \RSref{SecTP}.
\biglf
If the homogeneous
problem has a Maximum Principle, the well-posedness follows from it
via
$$
\|u\|_{\infty,\Omega}\leq
\|u\|_{\infty,\Gamma}=\|D(u)\|_{\infty,\Gamma}
\fa u\in U
$$
with the well-posedness norm $\|.\|_{WP}=\|.\|_{\infty,\Omega}$.
Trefftz methods for problems without a Maximum Principle need a
different way of proving well-posedness.
\biglf
The Method of Particular Solutions (MPS \RScite{atkinson:1985-1,
  zhu:1993-1,chen-et-al:1998-3,chen-et-al:1999-3, golberg:1996-1,
  muleshkov-et-al:2000-1,chen-et-al:2012-1}) will be shown
in Section \RSref{SecMPS} to inherit its well--posedness from
the well--posedness of the PDE problem.
\section{Approximation Problems}\RSlabel{SecAT}
We now forget operator equations until Section \RSref{SecTP}
and consider approximation problems
\eref{eqVapp0} on data spaces $F$.
These finite-dimensional linear approximation problems
clearly have solutions, but we are interested in the error
$\eta(f, F_r, F)$ in terms of the arguments.
In many cases, Approximation
Theory has good and handy results, but other situations may
be still open, e.g. the approximation by traces of
Fundamental Solutions, see Section \RSref{SecTP}. In general,
errors decrease with trial spaces getting larger and $f$ getting smoother,
at certain rates that are found in the literature.
\biglf
In general, users should try to get as much information on $u^*$
and $f=D(u^*)$ as possible, and then select trial spaces
$F_r=D(U_r)$ that approximate $f=D(u^*)$ well. It will be shown below
that the attainable approximation error dominates the error
in the operator equation solution, if stability issues
are handled properly. Remember that, in contrast to standard $h$-type finite
elements, the approach from Approximation Theory is free to choose
good trial spaces, and this freedom should be used wisely
and not be overdone.
\biglf
For an extreme case, consider papers
concerned with solving some PDE problem in 2D,
and providing an example with a true solution $u^*$
like
$$
u^*(x,y)=\exp(ax+by)
$$
for certain constants $a$ and $b$, or other cases where
one takes the exponential function of a low-degree polynomial.
There are plenty of such papers , e.g.
\RScite{kansa:2015-1,karageorghis-et-al:2017-1} for two
randomly chosen recent instances.
The solution can be approximated by a polynomial of low degree to
machine precision, the convergence being exponential in terms
of the degree, for the
function itself and all derivatives.
No matter what the PDE problem in the background is,
as long as polynomials are taken as trial functions,
one has a numerical problem of just a few degrees of freedom,
and one should not need many data
to get a solution up to very good accuracy, whatever the PDE problem is.
If there is no exponential convergence, something has gone wrong.
The same holds if the trial functions themselves have very good approximations
by polynomials, e.g. multiquadrics. Then one approximates
an approximate polynomial by approximate polynomials,
and this can only fail if the authors commit
certain numerical crimes that we have to consider in what follows.
It does not make sense to take trial spaces that are
not close to polynomials in such a case, but if non-polynomial trial spaces are
used, papers should contain an  experimental comparison to polynomials
that may outperform the trial space actually used. 
\section{Discretizing Approximation Problems}\RSlabel{SecStab}
We now reconsider approximation problems
\eref{eqVapp0} on data spaces $F$, but from a numerical perspective.
In view of Corollary \RSref{TheAPostErr}, we want a numerical method
that produces a function $\tilde f_r\in F_r$ with
\bql{eqffrCA}
\|f-\tilde f_r\|_F\leq C_A \|f-f_r^*\|_F
\eq
with a factor $C_A\geq 1$ that should be independent of $F_r$.
If this works, Corollary \RSref{TheAPostErr} yields the error bound
$$
\|u^*-\tilde u_r\|_{WP}\leq C_{WP}C_A\|f-f_r^*\|_F= C_{WP}C_A\eta(f,F_r,F)
$$
for $\tilde u_r\in U$ with $D(\tilde u_r)=\tilde f_r$ in terms
of the error provided by Approximation Theory, and we are done.
\biglf
But the problem with \eref{eqffrCA}
is that one cannot work directly on functions in Numerical Analysis.
The standard discretization in Numerical Analysis replaces functions by
finitely many of their values,
using a surjective {\em test} or
{\em sampling}  map $T_s\;:\;F\to V_s$ that
takes each function $f\in F$ into a vector
$T_s(f)$ in a {\em value space} $V_s$
of finitely many real numbers, e.g. values at points, or
integrals against test functions.
Of course, there is no stable recovery of $f$ from $T_s(f)$
because we have only finitely many data.
\biglf
Then \eref{eqVapp0} is replaced
by a {\em discrete} approximation problem
\bql{eqTsVapp} 
\displaystyle{\inf_{f_r\in F_r}\|T_s(f)-T_s(f_r)\|_{V_s}
  =\|T_s(f)-T_s(f_{r,s}^*)\|_{V_s}=:\eta(T_s(f),T_s(F_r),V_s)}
\eq
with a solution $f_{r,s}^*\in F_r$ that can actually be calculated
up to roundoff effects and numerical instabilities like bad choices of bases.
To allow some leeway, one may assume that one actually produces
a $\tilde f_{r,s}\in F_r$ with
\bql{eqtstildef}
\|T_s(f)-T_s(\tilde f_{r,s})\|_{V_s}\leq 2\eta(T_s(f),T_s(F_r),V_s).
\eq
When writing $\tilde f_{r,s}=D(\tilde u_{r,s})\in U_r$ due to $F_r=D(U_r)$,
the error bound \eref{eqapperrbnd} turns into
$$
  \|u^*-\tilde u_{r,s}\|_{WP}
  \leq  C_{WP}\|f-D(\tilde u_{r,s})\|_F
$$
  but we have no grip on $\|f-D(\tilde u_{r,s})\|_F$.
  Instead, we have $\|T_s(f)-T_s(D(\tilde u_{r,s}))\|_{V_s}$, but this is only a
  discrete norm. We need the transition from $\|T_s(f)\|_{V_s}$ to
  $\|f\|_{F}$, but this can only work on finite-dimensional subspaces. 
\biglf
Fortunately, Approximation Theory \RScite{rieger-et-al:2010-1}
helps with this again, because 
one can ask for 
a {\em stability inequality}
\bql{eqteststab}
\|f_r\|_F\leq C_{r,s}\|T_s(f_r)\|_{V_s} \fa f_r\in F_r
\eq
that inverts the sampling on the values of the trial space.
This allows to let the transition from the
full approximation problem \eref{eqVapp0} to the
discrete problem \eref{eqTsVapp}  
be stable in the sense
$$
\begin{array}{rcl}
  \|f- \tilde f_{r,s}\|_F
  &\leq &
   \|f- f_r^*\|_F+ \|f_r^*- \tilde f_{r,s}\|_F\\
  &\leq &
   \eta(f,F_r,F)+  C_{r,s}\|T_s(f_r^*)- T_s(\tilde f_{r,s})\|_{V_s}\\
 &\leq &
   \eta(f,F_r,F)+  C_{r,s}\|T_s(f_r^*)- T_s(f))\|_{V_s}
   +  C_{r,s}\|T_s(f)- T_s(\tilde f_{r,s})\|_{V_s}\\
 &\leq &
   \eta(f,F_r,F)+  3C_{r,s}\eta(T_s(f),T_s(F_r),V_s).
\end{array}
$$
\begin{theorem}\RSlabel{Thefull}
  Let a well-posed operator equation $D(u)=f$
  in the sense of Section \RSref{SecPDE}
  be given, and assume that for a trial space $U_r$ there is a test
  discretization on $F$ via sampling
  maps $T_s$ with \eref{eqteststab} for $F_r=D(U_r)$.
  Then the approximate solution $\tilde f_{r,s}=D(\tilde u_{r,s})$
  of the discretized approximation problem
  \eref{eqTsVapp} with \eref{eqtstildef} satisfies the error bound
  \bql{eqDiscOpErr}
  \|u^*-\tilde u_{r,s}\|_{WP}\leq C_{WP}
  \left(  \eta(f,F_r,F)+  3C_{r,s}\eta(T_s(f),T_s(F_r),V_s)   \right)
  \eq
  in terms of approximation errors. \qed
\end{theorem} 
This boils the problem down to stability inequalities
\eref{eqteststab} where we hope to bound $C_{r,s}$ independent of $r$ and $s$.
Note that this form of stability is necessary whenever one works
on finite values instead of functions and wants to conclude that
small discrete errors lead to small errors in function space.
The latter cannot be bypassed, because discrete norms will not
work in well-posedness inequalities. Any technique that
goes down to a finite system of equations or a finite
approximation problem will have to cope with such a stability argument,
but experience shows that authors only rarely care for the problem.
If the data error $\|T_sf-T_sD(\tilde u_{r,s})\|_{V_s}$ in an application paper
does not behave like what Approximation Theory predicts for 
$\|f-T_sD(u_r^*)\|_F=\eta(f, F_r,F)$, either $C_A$ or $C_{r,s}$
is not kept at bay, i.e. there is a flaw in the algorithm.
\section{Stability Inequalities}\RSlabel{SecSI}
These are an interesting and important part of Numerical Analysis,
and should be brought to the attention of a wider audience. 
We start with a seemingly simple classical case. 
\biglf
Assume a user wants to work in the space $C[-1,+1]$ 
with polynomials of order $M$, i.e. degree $M-1$.
All functions are replaced by their values on a set
$X_N:=\{x_1,\ldots,x_N\}\subset[-1,+1]$, and to let this identify
polynomials of order $M$ properly, one should let $X_N$ consist of $N\geq M$
different points, by the Fundamental Theorem of Algebra.
A {\em stability inequality} like
\eref{eqteststab} then is
$$
\|p\|_{\infty,[-1,+1]} \leq C_{M,X_N}\|p\|_{\infty,X_N}
$$
for all polynomials of order $M$, but what is the minimal stability
constant $C_{M,X_N}$?
\biglf
The answers of Approximation Theory are disppointing at first sight:
\begin{enumerate}
\item For $M=N$ and equidistant points,  $C_{M,X_N}$
  grows exponentially with $M$.
\item For $M=N$ and Chebyshev-distributed points,  $C_{M,X_N}$
  grows logarithmically with $M$.
\end{enumerate}
But if {\em oversampling} is used, the situation is {\em much} better: 
\begin{enumerate}
\item For $N\geq \pi M$ and Chebyshev-distributed points,  $C_{M,X_N}$
  is bounded independent of $M$. 
\item For $N \geq C M^2$ and equidistant points,  $C_{M,X_N}$
  is bounded independent of $M$.
\end{enumerate} 
The upshot is that replacing functions from an $M$--dimensional trial space
by $N\geq M$ function values is unstable unless
values are taken at well-chosen points and a serious amount of oversampling
is applied. For users solving operator equations, using
exotic trial spaces on nontrivial domains, this fact has to be taken into
account, because unbounded stability constants $C_{r,s}$ spoil
the approximation error rates in \eref{eqDiscOpErr}. The same holds if
users insist on having square linear systems. If these arise
from discretizing functions, instability {\em must} be expected.
However, if users choose examples with extremely smooth
true solutions that lead to very good or even exponential
convergence, the effect is not observable. 
\biglf
But a sufficient amount of oversampling can lead to uniform stability
under certain circumstances. We state a special case
of a result of \RScite{schaback:2015-4} based on
the extremely useful {\em norming set} notion of
\RScite{jetter-et-al:1999-1}:
\begin{theorem}\RSlabel{TheSupStab}
  For the space $F=C(\Omega)$ of continuous functions
  on a compact set $\Omega$ under the sup norm and all finite-dimensional
  subspaces $F_r$ there is a set $X_{s(r)}$ of points of $\Omega$ such that
  $C_{r,s(r)}\leq 2$. \qed
\end{theorem} 
A weak variant is
\begin{theorem}\RSlabel{TheL2Stab}
  For the space $F=L_2(\Omega)$ of square-integrable functions
  on a compact set $\Omega$ under the $L_2$ norm and all finite-dimensional
  subspaces $F_r$ there is a set of normalized test functionals
  defining a sampling operator $T_s$ 
  such that
  $C_{r,s(r)}\leq 2$. \qed
\end{theorem} 
Like in Section \RSref{SecPDE}, the difference
between the strong and weak case lies in what {\em data} means.
In Theorem \RSref{TheL2Stab}, normalized test functionals
are $L_2$ integrals against compactly supported test functions with norm 1.
\biglf
The above cases apply whenever the data space $F$ consists
exclusively of parts that
behave like $C(\Omega)$ or $L_2(\Omega)$. This fails if well-posedness is
stated in Hölder norms, but there is a bypass that will be treated elsewhere. 
\section{Trefftz Problems}\RSlabel{SecTP}
As stated in section \RSref{SecPDE}, a Trefftz problem for
a homogeneous differential equation with a Maximum Principle
reduces to approximation on $C(\Gamma)$ in the sup norm.
Thus Theorem \RSref{TheSupStab} applies to the stability problem,
implying
\begin{theorem}\RSlabel{TheTrefftz}
  Assume a well-posed problem for a homogeneous differential equation
  and Dirichlet boundary values, with the Maximum Principle being satisfied.
  Then for any trial space of homogeneous solutions there is a
  sufficiently fine
  set of test
  points on the boundary that guarantees uniform stability.\qed
  \end{theorem} 
It remains to check the approximation problem \eref{eqVapp0} in
$F=C(\Gamma)$. If played back to a trial space $U_r\subset U$
of homogeneous solutions, it turns into
$$
\displaystyle{\inf_{u_r\in U_r}\|u^*-u_r\|_{\infty,\Gamma}=
  :\eta(u^*_{|_\Gamma},{U_r}_{|_\Gamma},C(\Gamma))}
$$
and now the ball lies in the field of Approximation Theory, but the latter has
not much to say about this, unfortunately.
\biglf
If specialized to the MFS, the classical trial space consists of fundamental
solutions centered on a {\em fictitious boundary} outside of the domain,
but the approximation error is measured on the true boundary. This is a nasty
approximation problem that should get much more attention by
Approximation Theorists.
The papers
\RScite{schaback:2009-8,hon-schaback:2013-1} use special kernel-based trial
spaces where these approximation errors can be calculated, without any
fictitious boundary.
For the special case of equidistant points on
concentric circles and conformal images of
such configurations,
results of Katsurada 
\RScite{katsurada:1990-1,katsurada:1994-1},
handle the problem nicely  by Fourier analysis.
However, a good general theory is still missing.

\section{Method of Particular Solutions}\RSlabel{SecMPS}
Here one only has a differential operator as the data map $D\;:\;U\to F$,
and one works with pairs $(u_j,f_j)=(u_j,D(u_j))$ of trial functions
spanning trial spaces $U_r$ and $F_r=D(U_r)$, respectively. Then, given a
function $f\in F$, the approximation problem \eref{eqVapp0}
is posed, and this is completely independent of PDEs.
If an approximation $\tilde f_r\in F_r$ is found, one has
a function
$\tilde u_r\in U_r$ with
$D(\tilde u_r)=\tilde f_r$ that is taken as the desired result.
\biglf
However, the approximation problem \eref{eqVapp0} needs a discretization.
If carried out in $C(\Omega)$ with the sup norm,
we can invoke Theorem \RSref{TheSupStab},
implying
\begin{theorem}\RSlabel{TheMPS}
  Assume a differential operator
  $D\;:\;C^p(\Omega) \to F=C(\Omega)$ of order $p$ on a compact domain
  $\Omega\subset\R^d$. Then for any choice $(U_r,F_r)$ of
  trial spaces of particular solutions with $F_r=D(U_r)$,
  there is a finite set $X\subset\Omega$ such that
  the discretized approximation problem in the sup norm
  on $X$ is uniformly stable in the sense of Section \RSref{SecStab},
  i.e. the full approximation error is at most twice the
  discrete approximation error. \qed
\end{theorem} 
Often,
the choice of the $f_j$ is done first, in order to
use results on the approximation error by these functions,
but then one has to calculate the $u_j$ in order to transfer the approximation
back to $U$. In other cases, based on the smoothness of $u^*$, one can
use functions $u_j$ that give good approximation errors including
higher derivatives, and then the calculation of the $f_j$ is easy. 
\biglf
To compare these two in the kernel-based situation
using the Whittle--Mat\'ern kernels
generating Sobolev spaces $W_2^m(\R^d)$, we take the uniformly elliptic 
differential operator $D=Id-\Delta$ that maps isometrically
from $W_2^m(\R^d)$ to $W_2^{m-2}(\R^d)$ and back, by Fourier transform theory.
Let us assume that the true solution $u^*$ is in $W_2^m(\R^d)$.
If we start from translates of the kernel
of $W_2^m(\R^d)$ and form the images under $D$, we know
\RScite{wendland:2005-1} that the $L_\infty$ error in $F$ behaves like
$h^{m-2-d/2}$ where $h$ is the fill distance of the centers used to generate the
trial spaces. If we work backwards, we have to approximate $f=D(u^*)$,
but this is only in  $W_2^{m-2}(\R^d)$, and the best possible
approximation error in $L_\infty$ is again of order $h^{m-2-d/2}$.
This implies that finding the $u_j$ from the $f_j$ by complicated arguments
is likely
not to pay off.
The error is comparable in both cases. The MPS can be
effectively carried out from trial spaces $U_r$ in $U$, using
the spaces $F_r=D(U_r)$ for the approximation of $f$.
\section{Dual Reciprocity Method}\RSlabel{SecDRM}
But the Method of Particular Solutions ignores boundary conditions.
The standard application is a two-step technique for a problem
of the form $D=(L,B)$ with a differential operator $L$ and a
boundary operator
$B$, called {\em Dual Reciprocity Method}
\RScite{nardini-brebbia:1982-1,partridge-et-al:1992-1,chen-et-al:1999-1}.
\biglf
If the problem is posed as
$D(u)=(L(u),B(u))=f=(f_L,f_B)$, the MPS is applied first to come up with
an approximate solution of $L(u)=f_L$, i.e. with a function
$u_{MPS}$ such that
$$
\|f_L-L(u_{MPS})\|_F\leq \epsilon_{MPS}.
$$
The previous section dealt with this part, including stability
and error bounds.
\biglf
The second step takes the boundary values of $u_{MPS}$ and solves
the homogeneous problem $Lu=0$ with $B(u)=f_B-B(u_{MPS})$ by a Trefftz
or MFS technique. If we stabilize the approximation on the boundary
along the lines of Section \RSref{SecTP}, we get 
$$
\|B(\hat u)-B(\tilde u)\|_{\infty,\Gamma}\leq \epsilon_{T}
$$
for the approximate solution $\tilde u$ and the true solution $\hat u$
for the above homogeneous problem. 
\biglf
Note that both steps did not assume a full well-posedness.
If we now assume a well-posed Dirichlet problem
with the Maximum Principle
in the sense of \eref{eqinfWP},
our previous arguments imply 
$$
\begin{array}{rcl}
  \|u^*-\tilde u-u_{MPS}\|_{\infty,\Omega}
  &\leq &
  \|u^*-\tilde u-u_{MPS}\|_{\infty,\Gamma}
  +\|L(u^*-\tilde u-u_{MPS})\|_{\infty,\Omega}\\
  &=&
  \|B(\hat u)-B(\tilde u)\|_{\infty,\Gamma}
  +\|f_L-L(u_{MPS})\|_{\infty,\Omega}\\
  &\leq & \epsilon_{T}+\epsilon_{MPS}.
\end{array}
$$
\section{Direct Optimal Recovery}\RSlabel{SecDOR}
We viewed methods for
solving operator equations as an approximation of a function $u^*$
from their data $D(u^*)$. In section \RSref{SecPDE} we looked at
well--posed problems where $u^*$ is fully and stably determined by
the full data $D(u^*)$ comprising infinitely many conditions,
like in \eref{eqPoiStrong} and \eref{eqPoiWeak}. Later, in Section
\RSref{SecStab}, we went back to only partial and finite data
in order to have a numerically manageable problem.
\biglf
From an Approximation Theory viewpoint, this can be seen as a detour. 
If only finitely many data $\lambda_1(u^*),\ldots, \lambda_M(u^*)$
for linear functionals $\lambda_1,\ldots,\lambda_M\in U^*$ are given
right from the beginning,
e.g. a finite selection of the data functionals in
\eref{eqPoiStrong} or \eref{eqPoiWeak}, we should
find the best approximation to $u^*$ using this information only.
\biglf
Of course, this needs some regularization, and a simple way
\RScite{wu:1992-1} 
is to go into a suitable Hilbert space $H\subseteq U$
on which the functionals are
continuous and to construct the function $\tilde u\in H$ that has smallest
norm and shares the same data, i.e. satisfies
the generalized interpolation conditions
$\lambda_j(u^*)=\lambda_j(\tilde u),\;1\leq j\leq M$. By standard arguments,
the solution is a unique linear combination of the representer functions
$u_1,\ldots,u_M$ in $H$ of the functionals $\lambda_1,\ldots,\lambda_M$,
the coefficients being obtainable by solving a positive definite Gramian
matrix with entries $(\lambda_j,\lambda_k),\;1\leq j,k\leq M$.
\biglf
By another standard argument, the value $\tilde u(x)$ for any fixed
$x$ is the best linear prediction of any function value there, provided that
only the given data are available. Given $H$ and the data, there is no better
way of solving the problem pointwise. From a Machine Learning
viewpoint in Hilbert Spaces, this is an optimal way of learning
the solution of an operator equation from given training data.
\biglf
However, the method is not new at all. In the context of
kernel-based techniques, it is Symmetric Collocation
\RScite{wu:1992-1,
  fasshauer:1997-1,franke-schaback:1998-1,franke-schaback:1998-2a,
  schaback:2015-3}, but it can also be seen as a Rayleigh-Ritz method.
Due to its optimality properties, it is impossible to be outperformed
error--wise for the given data, 
but is has serious stability and complexity drawbacks that are hard to
overcome. A special case, connected to Trefftz methods and confined to
potential problems, is
in \RScite{schaback:2009-8,hon-schaback:2013-1}, but it deserves
extensions using new
kernels implementing singularity-free homogeneous solutions
of other differential operators. 
\bibliographystyle{plain}

\begin{thebibliography}{10}

\bibitem{atkinson:1985-1}
K.E. Atkinson.
\newblock The numerical evaluation of particular solutions for {P}oisson's
  equation.
\newblock {\em IMA Journal of Numerical Analysis}, 5:319--338, 1985.

\bibitem{Bogomolny:1985-1}
A.~Bogomolny.
\newblock Fundamental solutions method for elliptic boundary value problems.
\newblock {\em SIAM J. Numer. Anal.}, 22:644--669, 1985.

\bibitem{braess:2001-1}
D.~Braess.
\newblock {\em Finite Elements. Theory, Fast Solvers and Applications in Solid
  Mechanics}.
\newblock Cambridge University Press, 2001.
\newblock Second edition.

\bibitem{chen-et-al:2012-1}
C.~S. Chen, C.~M. Fan, and P.~H. Wen.
\newblock The method of approximate particular solutions for solving certain
  partial differential equations.
\newblock {\em Numer. Methods Partial Differential Equations}, 28(2):506--522,
  2012.

\bibitem{chen-et-al:1999-1}
C.S. Chen, C.A. Brebbia, and H.~Power.
\newblock Dual reciprocity method using compactly supported radial basis
  functions.
\newblock {\em Comm. Num. Meth. Eng.}, 15:137--150, 1999.

\bibitem{chen-et-al:2008-1}
C.S. Chen, A.~Karageorghis, and Y.S. Smyrlis.
\newblock {\em The Method of Fundamental Solutions - A Meshless Method}.
\newblock Dynamic Publishers, 2008.

\bibitem{chen-et-al:1999-3}
C.S. Chen, A.S. Muleshkov, and M.A. Golberg.
\newblock The numerical evaluation of particular solution for {P}oisson's
  equation - a revisit.
\newblock In C.A. Brebbia and H.~Power, editors, {\em Boundary Elements XXI},
  pages 313--322. WIT Press, 1999.

\bibitem{chen-et-al:1998-3}
C.S. Chen and Y.F. Rashed.
\newblock Evaluation of thin plate spline based particular solutions for
  {H}elmholtz-type operators for the {DRM}.
\newblock {\em Mech. Res. Comm.}, 25:195--201, 1998.

\bibitem{fairweather-et-al:1998-1}
G~Fairweather and A.~Karageorghis.
\newblock The method of fundamental solution for elliptic boundary value
  problems.
\newblock {\em Advances in Computatonal Mathematics}, 9:69--95, 1998.

\bibitem{fasshauer:1997-1}
G.~Fasshauer.
\newblock Solving partial differential equations by collocation with radial
  basis functions.
\newblock In A.~LeM\'ehaut\'e, C.~Rabut, and L.L. Schumaker, editors, {\em
  Surface Fitting and Multiresolution Methods}, pages 131--138. Vanderbilt
  University Press, Nashville, TN, 1997.

\bibitem{franke-schaback:1998-2a}
C.~Franke and R.~Schaback.
\newblock Convergence order estimates of meshless collocation methods using
  radial basis functions.
\newblock {\em Advances in Computational Mathematics}, 8:381--399, 1998.

\bibitem{franke-schaback:1998-1}
C.~Franke and R.~Schaback.
\newblock Solving partial differential equations by collocation using radial
  basis functions.
\newblock {\em Appl. Math. Comp.}, 93:73--82, 1998.

\bibitem{golberg:1996-1}
M.A. Golberg.
\newblock Recent developments in the numerical evaluation of particular
  solutions in the boundary element method.
\newblock {\em Appl. Math. Comp.}, 75:91--101, 1996.

\bibitem{golberg-et-al:1998-2}
M.A. Golberg and C.S. Chen.
\newblock The method of fundamental solutions for potential, {H}elmholtz and
  diffusion problems.
\newblock In M.A. Golberg, editor, {\em Boundary Integral Methods: Numerical
  and Mathematical Aspects}, pages 103--176. WIT Press, 1998.

\bibitem{hon-schaback:2013-1}
Y.C. Hon and R.~Schaback.
\newblock Solving the {3D} {L}aplace equation by meshless collocation via
  harmonic kernels.
\newblock {\em Adv. in Comp. Math.}, pages 1--19, 2013.

\bibitem{jetter-et-al:1999-1}
K.~Jetter, J.~St\"ockler, and J.D. Ward.
\newblock Error estimates for scattered data interpolation on spheres.
\newblock {\em Mathematics of Computation}, 68:733--747, 1999.

\bibitem{kansa:2015-1}
E.~Kansa.
\newblock Radial basis functions: achievements and challenges.
\newblock In {A.H.} Cheng and {C.} Brebbia, editors, {\em Boundary Elements and
  Other Mesh Reduction Methods XXXVII}, volume~61, pages 3--22. WIT
  Transactions on Modelling and Simulation, 2015.

\bibitem{karageorghis-et-al:2017-1}
A.~Karageorghis, M.~Jankowska, and C.S. Chen.
\newblock Kansa-{RBF} algorithms for elliptic problems in regular polygonal
  domains.
\newblock {\em Numerical Algorithms}, pages 1--23, Dec. 2017.

\bibitem{katsurada:1990-1}
M.~Katsurada.
\newblock Asymptotic error analysis of the charge simulation method in a jordan
  region with an analytic boundary.
\newblock {\em Journal of the Faculty of Science of Tokyo University, Section
  1A}, 37:635--657, 1990.

\bibitem{katsurada:1994-1}
M.~Katsurada.
\newblock Charge simulation method using exterior mapping functions.
\newblock {\em Japan Journal of Industrial and Applied Mathematics}, 11:47--61,
  1994.

\bibitem{kolodziej-zielinski:2009-1}
J.A. Kolodziej and A.P. Zielinski.
\newblock {\em Boundary Collocation Techniques and their Application in
  Engineering}.
\newblock WIT Press, 2009.

\bibitem{li-et-al:2008-1}
Z.-C. Li, T.-T. Lu, H.-Y. Hu, and A.~H.-D. Cheng.
\newblock {\em Trefftz and collocation methods}.
\newblock WIT Press, Southampton, 2008.

\bibitem{li-et-al:2010-1}
Z.C. Li, L.J. Young, H.T. Huang, and A.H.-D. Liu, Y.P.~Cheng.
\newblock Comparisons of fundamental solutions and particular solutions for
  {T}refftz methods.
\newblock {\em Eng. Anal. Bound. Elem.}, 34(3):248--258, 2010.

\bibitem{mathon-et-al:1977-1}
R.~Mathon and R.L. Johnston.
\newblock The approximate solution of elliptic boundary-value problems by
  fundamental solutions.
\newblock {\em SIAM J. Numer.\ Anal.}, 14:638--650, 1977.

\bibitem{muleshkov-et-al:2000-1}
A.S. Muleshkov, C.S. Chen, M.A. Golberg, and A.H-D. Cheng.
\newblock Analytic particular solutions for inhomogeneous {H}elmholtz-type
  equations.
\newblock In S.N. Atluri and F.W. Brust, editors, {\em Advances in
  Computational Engineering \& Sciences}, pages 27--32. Tech Science Press,
  2000.

\bibitem{nardini-brebbia:1982-1}
D.~Nardini and C.A. Brebbia.
\newblock A new approach to free vibration analysis using boundary elements.
\newblock In C.A. Brebbia, editor, {\em Boundary Element Methods in
  Engineering, Proc. 4th Int. Sem.}, pages 312--326. Springer-Verlag, 1982.

\bibitem{partridge-et-al:1992-1}
P.W. Partridge, C.A. Brebbia, and L.C. Wrobel.
\newblock {\em The Dual Reciprocity Boundary Element Method}.
\newblock CMP/Elsevier, 1992.

\bibitem{poullikkas-et-al:1998-1}
A.~Poullikkas, A.~Karageorghis, and G.~Georgiou.
\newblock The method of fundamental solutions for inhomogeneous elliptic
  problems.
\newblock {\em Computational Mechanics}, 22:100--107, 1998.

\bibitem{qin:2000-1}
Q.-H. Qin.
\newblock {\em The {T}refftz Finite and Boundary Element Method}.
\newblock WIT Press, 2000.

\bibitem{rieger-et-al:2010-1}
C.~Rieger, B.~Zwicknagl, and R.~Schaback.
\newblock Sampling and stability.
\newblock In M.~D\ae{}hlen, M.S. Floater, T.~Lyche, J.-L. Merrien,
  K.~M\o{}rken, and L.L. Schumaker, editors, {\em Mathematical Methods for
  Curves and Surfaces}, volume 5862 of {\em Lecture Notes in Computer Science},
  pages 347--369, 2010.

\bibitem{schaback:2008-7}
R.~Schaback.
\newblock An adaptive numerical solution of {MFS} systems.
\newblock In C.S. Chen, A.~Karageorghis, and Y.S. Smyrlis, editors, {\em The
  Method of Fundamental Solutions - A Meshless Method}, pages 1--27. Dynamic
  Publishers, 2008.

\bibitem{schaback:2009-8}
R.~Schaback.
\newblock Solving the {L}aplace equation by meshless collocation using harmonic
  kernels.
\newblock {\em Adv. in Comp. Math.}, 31:457--470, 2009.
\newblock DOI 10.1007/s10444-008-9078-3.

\bibitem{schaback:2010-2}
R.~Schaback.
\newblock Unsymmetric meshless methods for operator equations.
\newblock {\em Numerische Mathematik}, 114:629--651, 2010.

\bibitem{schaback:2015-4}
R.~Schaback.
\newblock All well--posed problems have uniformly stable and convergent
  discretizations.
\newblock {\em Numerische Mathematik}, 132:597--630, 2015.

\bibitem{schaback:2015-3}
R.~Schaback.
\newblock A computational tool for comparing all linear {PDE} solvers.
\newblock {\em Advances of Computational Mathematics}, 41:333--355, 2015.

\bibitem{wendland:2005-1}
H.~Wendland.
\newblock {\em Scattered Data Approximation}.
\newblock Cambridge University Press, 2005.

\bibitem{wu:1992-1}
Z.~Wu.
\newblock {H}ermite--{B}irkhoff interpolation of scattered data by radial basis
  functions.
\newblock {\em Approximation Theory and its Applications}, 8/2:1--10, 1992.

\bibitem{zhu:1993-1}
S.~Zhu.
\newblock Particular solutions associated with the {H}elmholtz operator used in
  {DRBEM}.
\newblock {\em Boundary Element Abstracts}, 4:231--233, 1993.

\end{thebibliography}


\end{document}